\documentclass{icm2010}
\usepackage [dvips]{graphicx}
\usepackage{psfrag}

\title[  ]{Conservative partially hyperbolic\\ dynamics
\thanks{Thanks to Christian Bonatti,  Keith Burns, Jordan Ellenberg, Andy Hammerlindl, Fran\c cois Ledrappier,  Charles Pugh, Mike Shub and Lee Wilkinson for reading earlier versions 
of this text and making several helpful suggestions.  This work was supported by the NSF.
}}
\author{Amie Wilkinson}

\contact[wilkinso@math.northwestern.edu]{Department of Mathematics\\ Northwestern University\\
2033 Sheridan Road\\ Evanston, IL 60208-2730, USA}

\theoremstyle{plain}
\newtheorem{main}{Theorem}
\newtheorem{coromain}{Corollary}
\newtheorem{conjecture}{Conjecture}

\newtheorem{theorem}{Theorem}[section]

\theoremstyle{remark}
\newtheorem{definition}[theorem]{Definition}

\newcommand{\RR}{{\mathbb R}}

\newcommand{\ZZ}{{\mathbb Z}}

\newcommand{\HH}{{\mathbb H}}

\newcommand{\Diff}{\operatorname{Diff}}

\newcommand{\PSL}{\operatorname{PSL}}
\newcommand{\torus}{{\mathbb T}}

\newcommand{\cF}{{\mathcal F}}

\newcommand{\cO}{{\mathcal O}}

\newcommand{\cU}{{\mathcal U}}

\newcommand{\cW}{{\mathcal W}}

\begin{document}

\begin{abstract} We discuss recent progress in understanding the dynamical  properties of partially hyperbolic diffeomorphisms that preserve volume.  The main topics addressed are density of stable ergodicity and stable accessibility, center Lyapunov exponents, pathological foliations, rigidity, and the surprising interrelationships between these notions.
\end{abstract}



\maketitle

\section*{Introduction}

Here is a story, told at least in part through the exploits of one of its main characters.  This character, like many a Hollywood (or Bollywood) star, has played a leading role in quite a few compelling tales;  this one ultimately concerns the dynamics of partially hyperbolic diffeomorphisms.

We begin with a connected, compact, smooth surface $S$ without boundary, of genus at least $2$.  The Gauss-Bonnet theorem tells us that the average curvature of any Riemannian metric on $S$ must be negative, equal to $2\pi\chi(S)$, where $\chi(S)$ is the Euler characteristic of $S$.  We restrict our attention to the metrics on $S$ of everywhere negative curvature; among such metrics, there is a finite-dimensional moduli space of {\em hyperbolic} metrics, which have constant curvature.  
Up to a normalization of the curvature, each hyperbolic surface may be represented  by a quotient $\HH/\Gamma$, where $\HH$ is the complex upper half plane with the  metric
$
y^{-2}(dx^2 + dy^2)
$,
 and
$\Gamma$ is a discrete subgroup of $\PSL(2,\RR)$, isomorphic to the fundamental group of $S$.
More generally, any negatively curved metric on $S$ lies in the conformal class of some hyperbolic metric, and the space of all such metrics is path connected.  Throughout this story, $S$ will be equipped with a negatively curved metric.

This negatively curved muse first caught the fancy of Jacques Hadamard in the late 1890's \cite{Ha898}.  Among other things, Hadamard studied the properties of geodesics on $S$ and a flow $\varphi_t\colon T^1S\to T^1 S$ on the unit tangent bundle to $S$ called the {\em geodesic flow}. The image  of a unit vector $v$ under the time-$t$ map of this flow is obtained by following the
unique unit-speed geodesic $\gamma_v\colon \RR\to S$ satisfying $\dot\gamma_v(0)=v$ for a distance $t$ 
and taking the tangent vector at that point: 
$$\varphi_t(v) := \dot\gamma_v(t).$$
This geodesic flow, together with  its close relatives, plays the starring role in the story told here.

\begin{figure}[h]
\begin{center}
\psfrag{f}{$\varphi$}
\includegraphics[scale=.5]{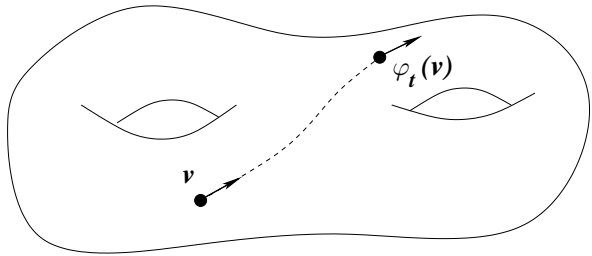}
\end{center}
\caption{The geodesic flow.}
\end{figure}

A theorem of Liouville implies that $\varphi_t$ preserves a natural probability measure $m$ on $T^1S$, known as {\em Liouville measure},
which locally is just the product of normalized area on $S$ with Lebesgue measure on the circle fibers.  Poincar{\'e} recurrence then  implies
that almost every orbit of the geodesic flow comes back close to itself infinitely often. 

In the special case where $S = \HH/\Gamma$ is a hyperbolic surface, the unit tangent bundle $T^1S$ is naturally
identified with $PSL(2,\RR)/\Gamma$, and the action of the geodesic flow $\varphi_t$ is realized by  left multiplication
by the diagonal matrix $$g_t = \left(\begin{array}{cc}e^{t/2}&0\\0 & e^{-t/2}\end{array}\right).$$
Liouville measure is normalized Haar measure.  

In his study of $\varphi_t$, Hadamard introduced the notion of the  {\em  stable manifold} of  a vector $v\in T^1S$:
$$
\cW^s(v) := \{w\in T^1S \mid  \lim_{t\to \infty} \text{dist}(\varphi_t(v),\varphi_t(w)) = 0\}.
$$
The proof that such  sets are manifolds is a nontrivial consequence of negative curvature and a noted accomplishment of Hadamard's.
Indeed, each stable manifold $\cW^s(v)$ is an injectively immersed, smooth copy of the real line, and taken together,
the stable manifolds form a foliation $\cW^s$ of $T^1M$. 
Similarly, one defines an {\em unstable manifold} by:
$$
\cW^u(v) := \{w\in T^1S \mid  \lim_{t\to -\infty}  \text{dist}(\varphi_t(v),\varphi_t(w)) = 0\}
$$ and denotes the corresponding unstable foliation $\cW^u$.
 The foliations $\cW^s$ and $\cW^u$ are key supporting players in this story.  

In the case where $S=\HH/\Gamma$, the stable manifolds are orbits of the {\em positive horocyclic flow} on 
$\PSL(2,\RR)/\Gamma$ defined by 
left-multiplication by 
 $$h^s_t = \left(\begin{array}{cc} 1 &t \\0 & 1\end{array}\right),$$ and the unstable
manifolds are orbits of the {\em negative horocyclic flow}, defined by
left-multiplication by 
$$h^u_t = \left(\begin{array}{cc} 1 &0 \\t & 1\end{array}\right).$$ 
This fact can be deduced from the explicit relations:
\begin{eqnarray}\label{e=matrixexpansion}
g_{-t} h^s_{r} g_{t} = h^s_{r e^{-t}} \quad\text{ and }\quad g_{-t} h^u_{r} g_{t} = h^u_{r e^{t}}.
\end{eqnarray}

The stable and unstable foliations stratify the future and past, respectively, of the geodesic flow. It might come as no surprise
that their features dictate the asymptotic behavior of the geodesic flow.  For example, Hadamard obtained from the existence of these foliations and Poincar{\'e} recurrence that periodic orbits for $\varphi_t$ are dense in $T^1S$.

Some 40 years after Hadamard received the Prix Poncelet for his work on surfaces, 
Eberhard Hopf introduced  a simple argument that
proved the ergodicity (with respect to Liouville measure)
of the geodesic flow on $T^1S$, for any closed negatively curved surface $S$ \cite{Ho39}. In particular, 
Hopf proved that almost every infinite geodesic in $S$ is dense (and uniformly distributed), not 
only in $S$, but in $T^1S$.  It was another thirty years before Hopf's result was extended
by Anosov to geodesic flows for negatively curved compact manifolds in arbitrary dimension.

Up to this point the discussion is quite well-known and classical, and from here the story can take many turns. For example, for arithmetic hyperbolic surfaces, the distribution of closed orbits of the flow and associated dynamical zeta functions quickly leads us into deep questions in analytic number theory.  Another path leads to the study the spectral theory of negatively curved surfaces, inverse problems and quantum unique ergodicity.  The path we shall take here  leads to the definition of partial hyperbolicity.

Let us fix a unit of time $t_0> 0$ and discretize the system $\varphi_t$ in these units; that is, we study
the dynamics of the time-$t_0$ map $\varphi_{t_0}$ of the geodesic flow.  From a digital age perspective  this is a natural thing to do; for example, to plot the orbits of a flow, a computer evaluates the flow at  discrete, usually equal, time intervals.  

If we carry this computer-based analogy one step further, we discover an interesting 
question. Namely, a computer does not ``evaluate the flow"  precisely, but rather uses an {\em approximation to the time-$t_0$ map} (such as an ODE solver or symplectic integrator) 
to compute its orbits.
To what extent does iterating this approximation retain the actual dynamical features of the flow $\varphi_t$, such as ergodicity?  

To formalize this question, we consider a {\em diffeomorphism} $f\colon T^1S\to T^1S$ such that the $C^1$ distance  $d_{C^1}(f,\varphi_{t_0})$ is small.  Note that $f$ in general will no longer embed in a flow.  While we assume that the distance from $f$ to $\varphi_{t_0}$ is small, this is no longer the case for the distance from $f^n$ to $\varphi_{n t_0}$, when $n$ is large.

\begin{figure}[h]
\begin{center}
\psfrag{x}{$x$}
\psfrag{f(x)}{$f(x)$}
\psfrag{phi(x)}{$\varphi_{t_0}(x)$}
\psfrag{fn(x)}{$f^n(x)$}
\psfrag{phit(x)}{$\varphi_{nt_0}(x)$}
\includegraphics[scale=1]{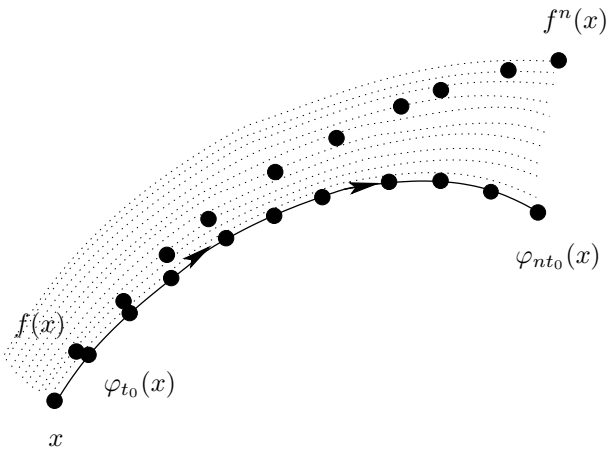}
\end{center}
\caption{$f^n(x)$ is not a good approximation to $\varphi_{nt_0}(x)$.}
\end{figure}

The earliest description of the dynamics of such a perturbation $f$ comes from a type of structural stability theorem proved by Hirsch, Pugh, and Shub~\cite{HPS77}.  The results there imply in particular that  if $d_{C^1}(f,\varphi_{t_0})$ is sufficiently small, then 
there exists an $f$-invariant \emph{center foliation} $\cW^c = \cW^c(f)$
that is homeomorphic to the orbit foliation $\cO$ of $\varphi_t$. The leaves of $\cW^c$ are smooth.  Moreover, the homeomorphism
$h\colon T^1S \to T^1S$ sending $\cW^c$ to $\cO$ is close to the identity and $\cW^c$ is the unique such foliation.

The rest of this paper is about $f$ and, in places, the foliation $\cW^c(f)$.  

What is known about $f$ is now substantial, but far from complete.  For example, the following basic problem is open.

\medskip

\noindent{\em Problem.} 
Determine whether $f$ has a dense orbit.  More precisely, does there exist a neighborhood $\cU$ of $\varphi_{t_0}$  in the space $\Diff^r(T^1S)$  of $C^r$ diffeomorphisms of $T^1S$ (for some $r\geq 1$) such that every $f\in \cU$ is topologically transitive?

\medskip

Note that $\varphi_{t_0}$ is ergodic with respect to volume $m$, and hence is itself topologically transitive.  In what follows, we will explain results from the last 15 years implying that any perturbation of $\varphi_{t_0}$ {\em that preserves volume} is ergodic, and hence has a dense orbit.  For perturbations that do not preserve volume,
 a seminal result of Bonatti and D\'iaz shows that $\varphi_{t_0}$ can be approximated arbitrarily well by $C^1$-open sets of transitive diffeomorphisms \cite{BD96}.  But the fundamental question of whether $\varphi_{t_0}$ lives in such an open set remains unanswered.

In most of the discussion here, we will work in  the conservative setting, in which the diffeomorphism $f$ preserves a volume probability measure.  To fix notation, $M$ will always denote a connected, compact Riemannian manifold without boundary, and $m$ will denote a probability volume on $M$.  For $r\geq 1$, we denote by $\Diff^r_m(M)$ the space of $C^r$ diffeomorphisms of $M$ preserving $m$, equipped with the $C^r$ topology.

\section{Partial hyperbolicity}

The map $\varphi_{t_0}$ and its perturbation $f$ are concrete examples of partially hyperbolic diffeomorphisms.  A diffeomorphism $f\colon M\to M$  of a compact Riemannian manifold $M$ is {\em partially hyperbolic} if there exists an integer $k\geq 1$ and a nontrivial, $Df$-invariant,  continuous splitting of the tangent bundle
$$
TM = E^s \oplus E^c \oplus E^u
$$
such that, for any $p\in M$ and unit vectors $v^s\in E^s(p)$,  $v^c\in E^c(p)$, and
 $v^u\in E^u(p)$:
$$\begin{array}{cccccc}
\|D_pf^k v^s\| & < & 1   &<& \|D_pf^k v^u\|, & \quad\text{and}  \\
\|D_pf^k v^s\|  & <  &\|D_pf^k v^c\| & < &\|D_pf^k v^u\|.  &
\end{array}
$$
Up to a change in the Riemannian metric, one can always take $k=1$ in this definition \cite{Go07}.
In the case where $E^c$ is the trivial bundle, the map $f$ is said to be {\em Anosov.} 
The central example $\varphi_{t_0}$ is partially hyperbolic: in that case, the bundle $E^c = \RR\dot\varphi$ is tangent to the orbits of the flow, and $E^s$ and $E^u$ are tangent to the leaves of $\cW^s$ and $\cW^u$, respectively.  

Partial hyperbolicity is a $C^1$-open condition: any diffeomorphism sufficiently
 $C^1$-close to a partially hyperbolic diffeomorphism is itself partially hyperbolic. Hence the
perturbations of $\varphi_{t_0}$ we consider are also partially hyperbolic. For an extensive discussion of examples of partially hyperbolic dynamical systems, see the survey articles \cite{BPSW, HP06,  HHU07s} and the book \cite{Pe04}. Among these examples are: the frame flow for a compact manifold of negative sectional curvature and most affine transformations of compact homogeneous spaces.

As is the case with the example $\varphi_{t_0}$, 
the stable and unstable bundles $E^s$ and $E^u$ of
an arbitrary partially hyperbolic diffeomorphism are always tangent to foliations, 
which we will again denote by $\cW^s$ and $\cW^u$ respectively; this is a consequence of partial hyperbolicity and a generalization of Hadamard's argument.
By contrast, the center bundle $E^c$ need not be tangent to a foliation, and can even be nowhere integrable.  In many  cases of interest, however, there is also a center 
foliation $\cW^c$ tangent to $E^c$: the content 
of the Hirsch-Pugh-Shub work in~\cite{HPS77} is the properties
of systems that admit such foliations, known as ``normally hyperbolic foliations."

There is a natural and slightly less general notion than integrability of $E^c$ that appears frequently in the literature.  We say that a partially hyperbolic diffeomorphism $f\colon M\to M$ is {\em dynamically coherent} if the subbundles $E^c\oplus E^s$ and  $E^c\oplus E^u$ are tangent to foliations $\cW^{cs}$ and $\cW^{cu}$, respectively,  of $M$.  If $f$ is dynamically coherent, then the center bundle $E^c$ is also integrable: one obtains the center foliation $\cW^c$ by intersecting the leaves of $\cW^{cs}$ and $\cW^{cu}$.  The examples $\varphi_{t_0}$ are dynamically coherent, as are their perturbations (by \cite{HPS77}: see \cite{BW08b} for a discussion).

\section{Stable ergodicity and the Pugh-Shub Conjectures}

Brin and Pesin \cite{BP74} and independently Pugh and Shub \cite{PS72}
first examined the ergodic properties of partially hyperbolic systems  in the early 1970's. 
The methods they developed give an ergodicity criterion for partially hyperbolic
$f\in \Diff^2_m(M)$ satisfying the following additional hypotheses:
\begin{itemize}
\item[(a)] the bundle $E^c$ is tangent to a $C^1$ foliation $\cW^c$, and 
\item[(b)] $f$ acts isometrically (or nearly isometrically) on the leaves of  $\cW^c$.
\end{itemize}
In \cite{BP74} it is shown that such an 
$f$ is ergodic with respect to $m$ if it satisfies a condition called accessibility.
\begin{definition}
A partially hyperbolic diffeomorphism $f:M\to M$ is {\em accessible} if
any point in $M$ can be reached from any other along an {\em $su$-path}, which
is a concatenation of finitely many subpaths, 
each of which lies entirely in a single leaf of $\cW^s$ or a single leaf of $\cW^u
$.
\end{definition}

This ergodicity criterion applies to the discretized geodesic flow $\varphi_{t_0}$: its center bundle
is tangent to the orbit foliation for $\varphi_t$, which is smooth, giving (a).  The action of $\varphi_{t_0}$ preserves the nonsingular vector field $\dot\varphi$, which implies (b).
It is straightforward to see that if $S$ is a hyperbolic surface, then $\varphi_{t_0}$ is accessible: the stable and unstable foliations are orbits of the smooth horocyclic flows
$h^s_t$ and $h^u_t$, respectively, and matrix multiplication on the level of the Lie algebra
$\frak{sl}_2$  shows that locally these flows
generate all directions in $\PSL(2,\RR)$:
\begin{eqnarray}\label{e=stabunst}
\frac12\left [\left(\begin{array}{cc}0 & 1\\ 0 & 0 \end{array}\right),  \left(\begin{array}{cc} 
0 & 0 \\ 1 & 0 \end{array}\right) \right ] =  \left(\begin{array}{cc} \frac12 & 0\\ 0 & -\frac12 \end{array}\right); 
\end{eqnarray}
the matrices appearing on the left are infinitesimal generators of the horocyclic flows, and the matrix on the right generates the geodesic flow.
Since $\varphi_{t_0}$ is accessible, it is ergodic. 

Now what of a small perturbation of $\varphi_{t_0}$?  As mentioned above, any $f\in \Diff_m^2(T^1S)$ sufficiently $C^1$ close to $\varphi_{t_0}$ also has a center foliation $\cW^c$, and the action of $f$ on the leaves is nearly isometric.  With some work, one can also show that $f$ is accessible (this was carried out in \cite{BP74}).  There is one serious reason why the ergodicity criterion of \cite{BP74} cannot be applied to $f$: the foliation $\cW^c$ is not $C^1$.  The leaves of $\cW^c$ are $C^1$, and the tangent spaces to the leaves vary continuously, but they do not vary smoothly. We will explore in later sections the extent to which  $\cW^c$ fails to be smooth, but for now suffice it to say that $\cW^c$ is pathologically bad, not only from a smooth perspective but also from a measure-theoretic one.

The extent to which $\cW^c$ is bad was not known at the time, but there was little hope of applying the existing techniques to perturbations of $\varphi_{t_0}$.
The first major breakthrough in understanding the ergodicity of perturbations of $\varphi_{t_0}$ came in the 1990's:
\begin{main}[Grayson-Pugh-Shub \cite{GPS94}] Let $S$ be a hyperbolic surface, and let $\varphi_t$ be the geodesic flow on $T^1S$.  Then $\varphi_{t_0}$ is {\em stably ergodic}:  there is a neighborhood $\cU$ of $\varphi_{t_0}$ in $\Diff_{m}^2(T^1 S )$ such that every $f\in\cU$ is ergodic  with respect to $m$.
\end{main}

The new technique introduced in \cite{GPS94} was a dynamical approach to understanding Lebesgue density points which they called {\em juliennes}.  The results in \cite{GPS94} were soon generalized to the case where $S$ has variable negative curvature \cite{Wi98} and to more general classes of partially hyperbolic diffeomorphisms \cite{PS97,PSh00}.
Not long after \cite{GPS94} appeared,  Pugh and Shub had formulated an influential circle of conjectures concerning the ergodicity of partially hyperbolic systems.

\begin{conjecture}[Pugh-Shub \cite{PS96}] \label{C=main}
On any compact manifold, ergodicity holds for an open and dense 
set of $C^2$ volume preserving partially hyperbolic diffeomorphisms.
\end{conjecture}

This conjecture can be split into two parts using the concept of accessibility.

\begin{conjecture}[Pugh-Shub \cite{PS96}] \label{C=main2}
Accessibility holds for an open and dense subset of $C^2$  
partially hyperbolic diffeomorphisms, volume preserving or not.
\end{conjecture}

\begin{conjecture}[Pugh-Shub \cite{PS96}] \label{C=main3}
A partially hyperbolic $C^2$ volume preserving diffeomorphism with 
the essential accessibility property is ergodic.
\end{conjecture}
 Essential accessibility is a measure-theoretic version
of accessibility that is implied by accessibility: $f$ is essentially accessible if for any two positive volume sets $A$ and $B$, there exists  an $su$-path in $M$ connecting some point in $A$ to some point in $B$ -- see \cite{BPSW} for a discussion.

In the next two sections, I will report on progress to date on these conjectures.

\medskip

\noindent{\em Further remarks.}
\begin{enumerate}
\item Volume-preserving Anosov diffeomorphisms (where $\dim E^c = 0$) are always ergodic.  This was proved by Anosov in his thesis \cite{An67}.  Note that Anosov diffeomorphisms are also accessible, since in that case the foliations $\cW^s$ and $\cW^u$ are transverse.  Hence all three conjectures hold true for Anosov diffeomorphisms.
\item It is natural to ask whether partial hyperbolicity is a necessary condition for stable ergodicity.  This is true when $M$ is $3$-dimensional \cite{DPU99} and also in the space of symplectomorphisms \cite{HT06, SX06}, but not in general \cite{Ta04a}.  What is true is that the related condition of having a dominated splitting is necessary for stable ergodicity (see \cite{DPU99}).
\item One can also ask whether for partially hyperbolic systems, 
 stable ergodicity implies accessibility.  
 If one works in a sufficiently high smoothness class, then this is not the case, as was shown in the groundbreaking paper of  F. Rodr\'iguez Hertz \cite{RH05}, who will also speak at this congress.  Hertz used methods from KAM theory to find an alternate route to stable ergodicity for certain essentially accessible systems.  
\item On the other hand, it is reasonable to expect that some form of accessibility is a necessary hypothesis for a general stable  ergodicity criterion for partially hyperbolic maps (see the discussion at the beginning of \cite{BW08b}).   Unlike Anosov diffeomorphisms, which are always ergodic, partially hyperbolic diffeomorphisms need not be ergodic.  For example, the product of an Anosov diffeomorphism with the identity map on any manifold is partially hyperbolic, but certainly not ergodic. See also  Theorem 11.16  in \cite{Beyond}.

\end{enumerate}

\section{Accessibility}

In general, the stable and unstable foliations of a partially hyperbolic diffeomorphism are not smooth (though they are not pathological, either -- see below).  Hence it is not possible in general  to use infinitesimal techniques to establish accessibility the way we did in equation (\ref{e=stabunst})  for the discretized hyperbolic geodesic flow.  The $C^1$ topology allows for enough flexibility in perturbations that Conjecture~\ref{C=main2} has been completely verified in this context:
\begin{main}[Dolgopyat-Wilkinson \cite{DW03}]\label{t=dw}  For any $r\geq 1$, accessibility holds for a $C^1$ open and  dense subset of the partially hyperbolic diffeomorphisms in $\Diff^r(M)$, volume-preserving or not.
\end{main}
Theorem~\ref{t=dw} also applies inside the space of partially hyperbolic symplectomorphisms.

More recently, the complete version of Conjecture~\ref{C=main3} has been verified for systems with $1$-dimensional center bundle.
\begin{main}[Rodr\'iguez Hertz-Rodr\'iguez Hertz-Ures \cite{HHU08b}]\label{t=rhrhu} For any $r\geq 1$, accessibility is $C^1$ open and $C^r$ 
dense among the partially hyperbolic diffeomorphisms in $\Diff^r_{m}(M)$ with one-dimensional center bundle.   
\end{main}
This theorem was proved earlier in a much more restricted context by Ni\c tic\u a-T\" or\" ok \cite{NT01}.  The $C^1$ openness of accesssibility was shown in \cite{Didier}. 
 A version of Theorem~\ref{t=rhrhu} for non-volume preserving diffeomorphisms was later  proved in \cite{BHHTU08}.  

The reason that it is possible to improve  Theorem~\ref{t=dw} from  $C^1$ density to $C^r$ density  in this context  is that the global structure of accessibility classes is well-understood.
By {\em accessibility class} we mean an equivalence class with respect to the relation generated by $su$-paths.  When the dimension of $E^c$ is $1$, accessibility classes are ($C^1$ immersed)
submanifolds.  Whether this is always true when $\dim(E^c)>1$ is unknown and is an important obstacle to attacking the general case of Conjecture~\ref{C=main2}.  

\medskip

\noindent{\em Further remarks.}
\begin{enumerate}
\item More precise criteria for accessibility have been established for special classes of partially hyperbolic systems such as discretized Anosov flows, skew products, and low-dimensional systems \cite{BW99, BPW00, HHU08a}.
\item Refined formulations of accessibility have been used  to study higher-order statistical properties of certain partially hyperbolic systems, in particular the discretized geodesic flow \cite{Do98, Li04}.  The precise relationship between accessibility and rate of mixing (in the absence of other hypotheses) remains a challenging problem to understand.
\item Accessibility also plays a key role in a recently developed Livsi\v c theory for partially hyperbolic diffeomorphisms,  whose conclusions closely mirror those in the Anosov setting \cite{KK96, Wliv}.
\end{enumerate}

\section{Ergodicity}\label{s=ergodicity}

Conjecture~\ref{C=main} has been verified under one additional, reasonably mild hypothesis:
\begin{main}[Burns-Wilkinson \cite{BW08a}]\label{t=bw} Let $f$ be $C^2$, volume-preserving,
partially hyperbolic and center bunched. 
If $f$ is essentially accessible, then $f$
is ergodic, and in fact has the Kolmogorov property.
\end{main}
The additional hypothesis is  ``center bunched."
A partially hyperbolic diffeomorphism 
$f$ is \emph{center bunched} if there exists an integer $k\geq 1$ such that for 
any $p\in M$ and any unit vectors $v^s\in E^s(p)$,  $v^c, w^c\in E^c(p)$, and  $v^u\in E^u(p)$:
\begin{eqnarray}\label{e=cb}
\|D_pf^k v^s\| \cdot \|D_p f^k w^c\|  <  \|D_p f^k v^c\| < \|D_p f^k v^u\| \cdot \|D_p f^k w^c\|.  
\end{eqnarray}
As with partial hyperbolicity, the definition of center bunching depends only on the smooth structure on $M$ and not the Riemannian structure;  if (\ref{e=cb}) holds for a given metric and $k\geq 1$, one can always find another metric for which (\ref{e=cb}) holds with $k=1$ \cite{Go07}.
In words, center bunching requires that the non-conformality of $Df\mid E^c$ be dominated by
the hyperbolicity of $Df \mid E^u\oplus E^s$.  Center bunching holds automatically if the
restriction of $Df$ to $E^c$ is conformal in some metric (for this metric, one can choose $k=1$).  In particular, if $E^c$ is one-dimensional, then $f$ is center bunched.  In the context where $\dim(E^c)=1$, Theorem~\ref{t=bw} was also shown in \cite{HHU08b}.

Combining Theorems~ \ref{t=rhrhu} and \ref{t=bw} we obtain:
\begin{coromain}\label{c=1donj}
The Pugh-Shub conjectures hold true among the partially hyperbolic diffeomorphisms with $1$-dimensional center bundle.
\end{coromain}
\medskip

\noindent{\em Further remarks.}
\begin{enumerate}
\item  The proof of Theorem~\ref{t=bw} builds on the original argument of Hopf for ergodicity of geodesic flows and incorporates a refined theory of the juliennes originally introduced in \cite{GPS94}.
\item It appears that the center bunching hypothesis in Theorem~\ref{t=bw} cannot be removed without a significantly new approach.  On the other hand, it is possible that Conjecture~\ref{C=main} will yield to other methods.
\item  Formulations of Conjecture~\ref{C=main} in the $C^1$ topology have been proved 
for low-dimensional center bundle \cite{BMVW04, RH05} and for symplectomorphisms \cite{ABW09}.  These formulations state that ergodicity holds for a {\em residual} subset in the $C^1$ topology.
\end{enumerate}

\section{Exponents}

By definition, a partially hyperbolic diffeomorphism produces uniform contraction and expansion in the directions tangent to $E^s$ and $E^u$, respectively.  In none of the results stated so far do we make any precise assumption on the growth of vectors in $E^c$ beyond the coarse bounds that come from partial hyperbolicity and center bunching.  In particular, an ergodic diffeomorphism in Theorem~\ref{t=bw} can have periodic points of different indices, corresponding to places in $M$ where $E^c$ is uniformly expanded, contracted, or neither.  The power of the julienne-based theory is that the hyperbolicity in $E^u\oplus E^s$, when combined with center bunching and accessibility, is enough to cause substantial mixing in the system, regardless of the precise features of the dynamics on $E^c$.

On the other hand, the asymptotic expansion/contraction rates in $E^c$ can give additional information about the dynamics of the diffeomorphism, and is a potentially important tool for understanding
partially hyperbolic diffeomorphisms that are not center bunched.

A real number $\lambda$ is a \emph{center Lyapunov exponent} of the partially
hyperbolic diffeomorphism $f: M\to M$ if there exists a nonzero vector
$v\in E^c$ such that
\begin{equation}\label{e=lyaplim}
\limsup_{n\to\infty} \frac{1}{n} \log \|Df^n(v)\| = \lambda.
\end{equation}
If $f$ preserves $m$, then Oseledec's theorem implies that the limit in \eqref{e=lyaplim}
exists for each $v\in E^c(x)$, for $m$-almost every $x$.
When the dimension of $E^c$ is $1$, the limit in \eqref{e=lyaplim} depends only on $x$,
and if in addition $f$ is ergodic with respect to $m$, then the limit takes a single
value $m$-almost everywhere. 

\begin{main}[Shub-Wilkinson \cite{SW00}]\label{t=SW}  There is an open set $\cU\subset \Diff^\infty_m(\torus^3) $ of  partially hyperbolic, dynamically coherent diffeomorphisms of the $3$-torus  $\torus^3 = \RR^3/\ZZ^3$ for which:
\begin{itemize}
\item the elements of $\cU$ approximate arbitrarily well (in the $C^\infty$ topology) the linear automorphism of $\torus^3$ induced by the matrix:
$$A = \left(\begin{array}{ccc} 2 & 1 & 0 \\1 & 1 & 0 \\ 0 & 0 & 1\end{array}\right)$$
\item the elements of $\cU$ are ergodic and  have positive center exponents, $m$-almost everywhere.
\end{itemize}
\end{main}

Note that the original automorphism $A$ has vanishing center exponents, {\em everywhere} on $\torus^3$, since $A$ is the identity map on the third factor.  Yet Theorem~\ref{t=SW} says that a small perturbation mixing the unstable and center directions of $A$ creates expansion in the center direction, almost everywhere on $\torus^3$.

The systems in $\cU$ enjoy the feature of being {\em non-uniformly hyperbolic}: the Lyapunov exponents in every direction (not just center ones) are nonzero, $m$-almost everywhere.  The  well-developed machinery  of Pesin theory guarantees a certain level of chaotic behavior from nonuniform hyperbolicity alone.  For example, a nonuniformly hyperbolic diffeomorphism has at most countably many ergodic components, and a mixing partially hyperbolic diffeomorphism is Bernoulli (i.e. abstractly isomorphic to a Bernoulli process).  A corollary of Theorem~\ref{t=SW} is that the elements of $\cU$ are Bernoulli systems.

The constructions in \cite{SW00} raise the question of whether it might be possible to ``remove zero exponents" from any partially hyperbolic diffeomorphism via a small perturbation.  If so, then one might be able to bypass the julienne based theory entirely and use techniques from Pesin theory instead as an approach to Conjecture~\ref{C=main}.
More generally, and wildly optimistically, one might ask whether {\em any} $f\in \Diff_m^2(M)$ with at least one nonzero Lyapunov exponent on a positive measure set might be perturbed to produce nonuniform hyperbolicity on a positive measure set (such possibilities are discussed in \cite{SW00}).

There is a partial answer to these questions in the $C^1$ topology, due to Baraveira and Bonatti \cite{BB03}.  The results there imply in particular that if $f \in\Diff^r_m(M)$ is partially hyperbolic, then there exists $g \in\Diff^r_m(M)$, $C^1$-close to $f$ so that the {\em sum} of the center Lyapunov exponents is nonzero.

\medskip

\noindent{\em Further remarks.}
\begin{enumerate}
\item Dolgopyat proved that the same type of construction as in \cite{SW00}
can be applied to the discretized geodesic flow $\varphi_{t_0}$ for a negatively curved surface $S$ to produce perturbations with nonzero center exponents \cite{Dol04}.  See also \cite{Ru03}
for further generalizations of \cite{SW00}.
\item  An alternate approach to proving Conjecture~\ref{C=main} has been proposed, taking into account the center Lyapunov exponents \cite{BDP02}.
For systems with $\dim(E^c)=2$, this program has very recently been carried out in the $C^1$ topology in \cite{HHTU10}, using a novel application of the technique of blenders, a concept introduced in \cite{BD96}.
\end{enumerate}

\section{Pathology}

There is a curious by-product of nonvanishing Lyapunov exponents for the open set $\cU$ of examples in Theorem~\ref{t=SW}.  
By \cite{HPS77}, there is a center foliation $\cW^c$ for each $f\in \cU$,
homeomorphic to the trivial $\RR/\ZZ$ fibration of $\torus^3 = \torus^2\times \RR/\ZZ$;
in particular, the center leaves are all compact.  The almost everywhere exponential growth
associated with nonzero center exponents is incompatible with the compactness of the
center foliation, and so the full volume set with positive center exponent must meet
almost every leaf in a zero set (in fact a finite set \cite{RW01}).

The same type of phenomenon occurs in perturbations of the discretized geodesic flow $\varphi_{t_0}$.  While in that case the leaves of $\cW^c$ are mostly noncompact, they are in a sense ``dynamically compact." An adaptation of the arguments in  \cite{RW01} shows that any
perturbation of $\varphi_{t_0}$ with nonvanishing center exponents, such as those constructed by Dolgopyat in \cite{Dol04}, have atomic disintegration of volume along center leaves.
\begin{definition}
A foliation $\cF$ of $M$ with smooth leaves has {\em atomic disintegration of volume} along its leaves if there exists $A\subset M$ such that
\begin{itemize}
\item $m(M\setminus A) = 0$, and
\item $A$ meets each leaf of $\cF$ in a discrete set of points (in the leaf topology).
\end{itemize}
\end{definition}

\begin{figure}[h]
\begin{center}
\leavevmode
\includegraphics[scale=.35]{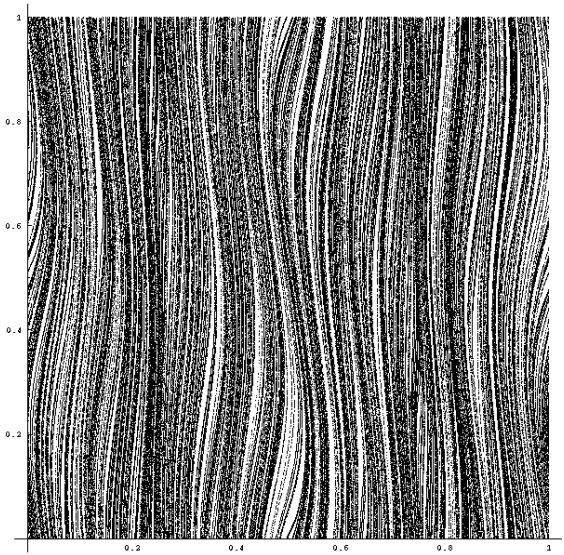}
\end{center}
\caption{A pathological foliation}
\end{figure}

At the opposite end of the spectrum from atomic disintegration of volume is a property called absolute continuity.
A foliation  $\cF$ is {\em  absolutely continuous} if holonomy maps between
smooth transversals send zero volume sets to zero volume sets.   If $\cF$
has smooth leaves and is absolutely continuous, then for every set $A\subset M$
satisfying  $m(M\setminus A) = 0$, the intersection of $A$ with the leaf $\cF$
through $m$-almost every point in $M$ has full leafwise Riemannian volume. In this sense
Fubini's theorem holds for absolutely continuous foliations.  If $\cF$ is a $C^1$ foliation, then it is absolutely continuous, but absolute continuity is a strictly weaker property.

Absolute continuity has long played a central role in
smooth ergodic theory.  Anosov and Sinai~\cite{An67,AS67} proved in the 60's  that the stable and unstable foliations of globally hyperbolic (or Anosov) systems are absolutely continuous,
even though they fail to be $C^1$ in general. Absolute continuity was a key ingredient in
Anosov's celebrated proof~\cite{An67} that the geodesic flow for any compact, negatively curved manifold is ergodic.  When the center foliation for $f$ fails to be absolutely continuous,
this  means that one cannot ``quotient out by the center direction"
to study ergodic properties $f$.

The existence of such pathological center foliations was first demonstrated by A. Katok (whose construction was written up by Milnor in \cite{Mi97a}).  Theorem~\ref{t=SW} shows that this type of pathology can occur in open sets of diffeomorphisms and so is inescapable in general.  In the next section, we discuss the extent to which this pathology is the norm.

\medskip

\noindent{\em Further remarks.}
\begin{enumerate}
\item  An unpublished letter of Ma\~n\'e to Shub examines the consequences of nonvanishing Lyapunov center exponents on the disintegration of volume along center foliations. Some of the ideas there are pursued in greater depth in \cite{HP07}.  
\item The examples of Katok in \cite{Mi97a} in fact have center exponents almost everywhere equal to 0, showing that nonvanishing center exponents is not a necessary condition for atomic disintegration of volume.  
\item  Systems for which the center leaves are not compact (or even dynamically compact) also exhibit non-absolutely continuous center foliations, but the disintegration appears to be potentially much more complicated than just atomic disintegration \cite{SX09, Go10}.
\end{enumerate}

\section{Rigidity}\label{s=rigidity}

Examining in greater depth the potential pathologies of center foliations, we discover a rigidity
phenomenon.  To be concrete, let us consider the case of a perturbation $f\in \Diff^\infty_m(M)$ of the discretized geodesic flow on a negatively-curved surface.  If the perturbation $f$ happens to be the time-one map of a smooth flow, then $\cW^c$ is the orbit foliation for that flow.  In this case the center foliation for $f$ is absolutely continuous -- in fact, $C^\infty$.
In general, however, a perturbation $f$ of $\varphi_{t_0}$ has no reason to embed in a smooth flow, and one can ask how the volume $m$ disintegrates along the leaves of $\cW^c$.

There is a complete answer to this question:
\begin{main}[Avila-Viana-Wilkinson \cite{AVW1}]\label{t=avw}
Let $S$ be a closed negatively curved surface, and let $\varphi_t: T^1S\to T^1S$ be the geodesic flow.

For each $t_0>0$, there is a neighborhood $\cU$ of $\varphi_{t_0}$ in  $\Diff^\infty_m(T^1S)$ such that for each $f\in \cU$:
\begin{enumerate}
\item either $m$ has atomic disintegration along the center foliation $\cW^c$, or
\item $f$ is the time-one map of a $C^\infty$, $m$-preserving flow.
\end{enumerate}
\end{main}
What Theorem~\ref{t=avw} says  is that,  in this context, nothing lies between $C^\infty$ and absolute singularity of $\cW^c$ -- {\em pathology is all that can happen}.  The geometric measure-theoretic properties of $\cW^c$ determine completely a key dynamical property of $f$ -- whether it embeds in a flow.

The heart of the proof of Theorem~\ref{t=avw}
is to understand what happens when the center Lyapunov exponents \emph{vanish}.
For this, we use tools that originate in the study of random matrix products.
The general theme of this work, summarized by Ledrappier in \cite{Le86} is that
``entropy is smaller than exponents, and entropy zero implies  deterministic."
Original results concerning the Lyapunov exponents of random matrix products,
due to Furstenberg, Kesten~\cite{FK60, Fu63}, Ledrappier~\cite{Le86}, and others,
have been extended in the past decade to deterministic products of linear cocycles
over hyperbolic systems by Bonatti, Gomez-Mont, Viana~\cite{BGV03,BoV04,Almost}.
The Bernoulli and Markov measures associated with random products in those earlier
works are replaced in the newer results by invariant measures for the hyperbolic
system carrying a suitable product structure.

Recent work of Avila, Viana~\cite{AV3} extends this hyperbolic theory from linear to
\emph{diffeomorphism} cocycles, and these results are used in a central way.
For cocycles over volume preserving partially hyperbolic systems, Avila, Santamaria,
and Viana \cite{ASV} have also recently produced related results, for both linear and
diffeomorphism cocycles, which also play an important role in the proof.  The proof
in \cite{ASV} employs julienne based techniques, generalizing the 
arguments in \cite{BW08b}.

\medskip

\noindent{\em Further remarks.}
\begin{enumerate}
\item  The only properties of $\varphi_{t_0}$ that are used in the proof of Theorem~\ref{t=avw} are accessibility, dynamical coherence, one-dimensionality of $E^c$, the fact that $\varphi_{t_0}$ fixes the leaves of $\cW^c$, and $3$-dimensionality of $M$. 
There are also more general formulations of Theorem~\ref{t=avw} in \cite{AVW1} that relax these hypotheses in various directions.  For example, a similar result holds for systems in dimension 
$3$ for whom all center manifolds are compact.
\item Deep connections between Lyapunov exponents and geometric properties of invariant measures have long been understood \cite{Le84a,LY85a,LY85b,Ka80,BPS99}.
Theorem~\ref{t=avw} establishes new connections in the partially hyperbolic context.
\item Theorem~\ref{t=avw} gives conditions under which one can recover the action of a Lie group (in this case $\RR$) from that of a discrete subgroup (in this case $\ZZ$).
These themes have arisen in the related context of measure-rigidity for algebraic partially
hyperbolic actions by Einsiedler, Katok, Lindenstrauss~\cite{EKL06}.
It would be interesting to understand more deeply the connections between these works.
\end{enumerate}

\section{Summary, questions.}

We leave this tale open-ended, with a few questions that have arisen naturally in its course.

\medskip

\noindent{\bf New criteria for ergodicity.}  Conjecture~\ref{C=main} remains open.  As discussed in Section~\ref{s=ergodicity}, the julienne based techniques using the Hopf argument might have reached their limits in this problem (at least this is the case in the absence of a significantly new idea).   One alternate approach which seems promising employs Lyapunov exponents and blenders \cite{HHTU10}.   Perhaps a new approach will find a satisfying conclusion to this part of the story.

\medskip

\noindent{\bf Classification problem.} A basic question is to understand which manifolds support partially hyperbolic diffeomorphisms.  As the problem remains open in the classical Anosov case (in which $E^c$  is zero-dimensional), it is surely exremely difficult in general.   There has been significant progress in dimension $3$, however; for example, using techniques in the theory of codimension-$1$ foliations, Burago and Ivanov proved that there are no partially hyperbolic diffeomorphisms of the $3$-sphere \cite{BI08}.

Modifying this question slightly, one can ask whether the partially hyperbolic diffeomorphisms in low dimension must belong to certain ``classes" (up to homotopy, for example) -- such as time-$t$ maps of flows, skew products, algebraic systems, and so on.  Pujals has proposed such a program in dimension $3$, which has spurred several papers on the subject \cite{BBI04, BW05, HHU08a, Hammerlindl}.

It is possible that if one adds the hypotheses of dynamical coherence and absolute continuity of the center foliation, then there is  such a classification.  Evidence in this direction can be found in \cite{AVW1}.

\medskip

\noindent{\bf Nonuniform and singular partial hyperbolicity.}  Unless all of its Lyapunov exponents vanish almost everywhere, {\em any} volume-preserving diffeomorphism is in some sense ``nonuniformly partially hyperbolic."  Clearly such a general class of systems will not yield to a single approach.  Nonetheless, the techniques developed recently are quite powerful and should shed some light on certain systems that are close to being partially hyperbolic.  Some extensions beyond the uniform setting have been explored in \cite{ABW09}, in which the center bunching hypotheses in \cite{BW08b} has been replaced by a pointwise, nonuniform center bunching condition.  This gives new classes of stably ergodic diffeomorphisms that are 
not center bunched.

 It is conceivable that the methods
in \cite{ABW09} may be further extended to apply in certain ``singular
partially hyperbolic'' contexts where
partial hyperbolicity holds on an open, noncompact subset of the manifold $M$
but decays in strength near the boundary.  Such conditions hold, for
example, for geodesic flows on certain nonpositively curved manifolds.
Under suitable accessibility hypotheses, these systems should be
ergodic with respect to volume.

\medskip

\noindent{\bf Rigidity of partially hyperbolic actions.}  The rigidity phenomenon described in Section~\ref{s=rigidity} has only begun to be understood. To phrase those results in a more general context, we consider a smooth, nonsingular action of an abelian Lie group $G$ on a manifold $M$.  Let $H$ be a discrete group acting on $M$, commuting with the action of $G$,
and whose elements are partially hyperbolic diffeomorphisms in $\Diff^\infty_m(M)$. 
Can such an action be perturbed, preserving the absolute continuity of the center foliation? 
How about the elements of the action?  When absolute continuity fails, what happens?

The role of accessibility and accessibility classes has been exploited in a serious way in the important work of  Damjanovi{\'c} and A. Katok on rigidity of abelian actions on quotients of $\operatorname{SL}(n,\RR)$  \cite{DaKa09b}.  It seems reasonable that these explorations can be pushed further, using some of the techniques mentioned here, to prove rigidity results for other partially hyperbolic actions.  A simple case  currently beyond the reach of existing methods is to understand perturbations  of the action of a $\ZZ^2$ lattice in the diagonal subgroup on $\operatorname{SL}(2,\RR) \times \operatorname{SL}(2,\RR)/\Gamma$, where $\Gamma$ is  an irreducible lattice.

 Our final question takes us further afield, but back once again to the geodesic flow.
Fix a closed hyperbolic surface $S$, and consider the standard action on $T^1S$  by the upper triangular  subgroup  $T < \PSL(2,\RR)$, which contains both the geodesic and positive horocyclic flows.   Ghys proved that this action is highly rigid and admits no $m$-preserving $C^\infty$ deformations \cite{Gh85}.  Does the same hold true for some countable subgroup of $T$?  For example, consider the solvable Baumslag Solitar subgroup $BS(1,2)$
 generated by the elements
$$a = \left(\begin{array}{cc} \sqrt{2} & 0 \\  0& \frac{1}{\sqrt{2}}\end{array}\right)
 \quad\text{ and }\, 
b= \left(\begin{array}{cc} 1 & 1 \\0 & 1\end{array}\right),
$$
which has the presentation $BS(1,2) = \langle a,b \,\mid\,  aba^{-1} = b^2 \rangle$. Can the standard action be perturbed inside of $\Diff^\infty_m(T^1S)$?
More generally, can one classify all faithful representations 
$$\rho\colon BS(1,2) \to \Diff^\infty_m(M),
$$
where $M$ is a $3$-manifold?  For results of a similar nature in lower dimensions, see
\cite{BursW04, Poltero02}.

\bibliographystyle{amsplain}
\providecommand{\bysame}{\leavevmode\hbox to3em{\hrulefill}\thinspace}
\providecommand{\MR}{\relax\ifhmode\unskip\space\fi MR }
\providecommand{\MRhref}[2]{%
  \href{http://www.ams.org/mathscinet-getitem?mr=#1}{#2}
}
\providecommand{\href}[2]{#2}

\end{document}